\documentclass[11pt]{article}
\newcommand{\documentdate}{4 IV 2021}

\usepackage{a4wide,latexsym,amsmath}

\title{H\" {o}lder Gradient Descent and Adaptive Regularization Methods in Banach Spaces for First-Order Points}

\author{
   S. Gratton%
   \thanks{Universit\'e de Toulouse, INP, IRIT, Toulouse, France. Email:
     serge.gratton@enseeiht.fr. Work partially supported by 3IA Artificial and
     Natural Intelligence Toulouse Institute, French "Investing for the Future
     - PIA3" program under the Grant agreement ANR-19-PI3A-0004"}, 
   ~S. Jerad%
   \thanks{ANITI, Universit\'e de Toulouse, INP, IRIT, Toulouse, France. Email:
     sadok.jerad@enseeiht.fr}
   ~and Ph. L. Toint%
   \thanks{NAXYS, University of Namur, Namur, Belgium. Email: philippe.toint@unamur.be}
}
\newcommand{\beqn}[1]{\begin{equation}\label{#1}}
\newcommand{\eeqn}{\end{equation}}
\newcommand{\req}[1]{(\ref{#1})}

\newcommand{\tim}[1]{\;\; \mbox{#1} \;\;}

\newtheorem{theorem}{Theorem}[section]
\newtheorem{lemma}[theorem]{Lemma}

\newcommand{\numsection}[1]{\section{#1}\setcounter{equation}{0}}

\newcounter{algo}[section]
\renewcommand{\thealgo}{\thesection.\arabic{algo}}
\newcommand{\llem}[2]{\vspace{\baselineskip} 
\noindent\framebox[\textwidth]{\parbox{0.95\textwidth}{
\begin{lemma} \label{#1} \rm #2 \end{lemma} } } \vspace{\baselineskip} }
\newcommand{\algo}[3]{\refstepcounter{algo}
\begin{center}\begin{figure}[htbp]
\framebox[\textwidth]{
\parbox{0.95\textwidth} {\vspace{\topsep}
{\bf Algorithm \thealgo : #2}\label{#1}\\
\vspace*{-\topsep} \mbox{ }\\
{#3} \vspace{\topsep} }}
\end{figure}\end{center}}
\newcommand{\bpr}{{\bf Proof.} \hspace{1.5mm}}
\newcommand{\epr}{\hfill $\Box$ \vspace*{1em}}
\newcommand{\proof}[1]{
\begin{list}{}{
\setlength{\topsep}{0.0pt}
\setlength{\partopsep}{0.0pt}
\setlength{\leftmargin}{0.025\textwidth}
\setlength{\rightmargin}{0.5\leftmargin}
\setlength{\labelwidth}{0.5\leftmargin}
\setlength{\labelsep}{0.25\leftmargin}}
\item \bpr #1 \epr \noindent
\end{list}}
\newcommand{\lthm}[2]{\vspace{\baselineskip} 
\noindent\framebox[\textwidth]{\parbox{0.95\textwidth}{
\begin{theorem} \label{#1} \rm #2 \end{theorem} } } \vspace{\baselineskip} }
\newcommand{\pd}[1]{\langle #1 \rangle}
\newcommand{\ii}[1]{\{ 1, \ldots, #1 \}}

\newcommand{\calD}{{\cal D}} 
\newcommand{\calF}{{\cal F}} 
\newcommand{\calO}{{\cal O}} 
\newcommand{\calS}{{\cal S}} 
 
\newcommand{\calV}{{\cal V}}
\newcommand{\calH}{\mathcal{H}}
\renewcommand{\Re}{\hbox{I\hskip -2pt R}}
\newcommand{\bigfrac}[2]{\frac{\displaystyle #1}{\displaystyle #2}}
\newcommand{\sfrac}[2]{{\scriptstyle \frac{#1}{#2}}}

\newcommand{\eqdef}{\stackrel{\rm def}{=}}
\newcommand{\bigsum}{\displaystyle \sum}
\newcommand{\al}[1]{{\footnotesize{\sf #1}}}
\newcommand{\tal}[1]{{\normalsize {\sf #1}}}
\newcommand{\np}[1]{\|#1\|_{\cal V}}
\newcommand{\nd}[1]{\|#1\|_{\cal V'}}
\date{\documentdate}

\begin{document}

\maketitle

\begin{abstract}
This paper considers optimization of smooth nonconvex functionals in smooth
infinite dimensional spaces. A H\"older gradient descent algorithm is first
proposed for finding approximate first-order points of regularized polynomial
functionals. This method is then applied to analyze the evaluation complexity of
an adaptive regularization method which searches for approximate first-order
points of functionals with $\beta$-H\"older continuous derivatives. It is shown
that finding an $\epsilon$-approximate first-order point requires at most
$\calO(\epsilon^{-\frac{p+\beta}{p+\beta-1}})$ 
evaluations of the functional and its first $p$ derivatives.
\end{abstract}


\numsection{Introduction}

The analysis of adaptive regularization (AR) algorithms for nonlinear (and
potentially nonconvex) optimization has been a very active field in recent
years (see
\cite{Grie81,NestPoly06,CartGoulToin09a,CartGoulToin11d,CartGoulToin12a,%
BianLiuzMoriScia15,GrapYuanYuan15a,BianScia16,BirgGardMartSantToin17,Mart17,%
GratRoyeVice17,BergDiouGrat17,BellGuriMoriToin19,CartGoulToin20b},
to cite only a few).  This sustained interest of the research
community is motivated in part by the fact that these methods not only work
well in practice, but also exhibit excellent worst-case evaluation complexity
bounds: one can indeed prove that the number of function and derivatives
evaluations which may be required to find an approximate critical point is
small, at least compared to similar bounds for other standard methods such as
linesearch-based Newton or trust-region algorithms
\cite{NestPoly06,CartGoulToin11d}.  As it turns out, evaluation complexity
results obtained for AR methods and nonconvex problems have been obtained, to
the best of the authors' knowledge, in the context of $\Re^n$.  It is the
purpose of this short note to show that this need not be the case, and that
evaluation complexity bounds for computing approximate first-order critical
point can be derived in infinite-dimensional Banach spaces.

The motivation for this generalization is a matter of coherence when
optimization algorithms are applied to large-scale discretized problems: 
it is then important to show that AR methods continue to make sense in the limit,
as the discretization mesh converges to zero. This coherence, sometimes called
``mesh independence'', has long been considered as an important feature of
numerical optimization methods
\cite{KellSach87,AllgBohmPotrRhei86,DeufPotr92,Hein93,UlbrUlbr00}.  For
trust-region methods, this was studied in \cite{Toin88} in the Hilbert space 
context, and developed for Hilbert and Banach spaces in
\cite[Section~8.3]{ConnGoulToin00}. Considering the question for AR algorithms
therefore seems a natural development in this line of research.

The outline of adaptive regularization methods is today quite well-known for
finite dimensional spaces (see \cite{BirgGardMartSantToin17}, for instance),
but difficulties arise in the nonconvex infinite dimensional space case.  The main
problem is that the existence of a suitable step at a given iteration of the
method typically hinges on approaching a minimizer of the regularized model,
which may no longer exist in infinite dimensions.  Our analysis circumvents
that problem by proposing a specialized optimization technique which guarantees
an acceptable step. 

\noindent
{\bf Contributions.} Having set the scene, we now make our contribution more
precise.
\begin{itemize}
\item We first analyse the convergence of a method for minimizing
  polynomial functionals with a general differentiable convex regularization
  whose gradients satisfy a generalized H\"older condition. To our knowledge, no
  such regularization has been considered before, even in finite dimensional spaces.
\item We then propose an adaptive regularization algorithm for finding
  first-order points of nonconvex functions having H\"older continuous $p$-th
  derivative (in the Fr\'{e}chet sense) and analyze its evaluation
  complexity. We show that the sharp complexity bound known
  \cite{CartGoulToin18a} for the finite-dimensional case is recovered, in
  that the algorithm requires at most $\calO\big(\epsilon^{-\frac{p+\beta}{p+\beta-1}}\big)$ 
  evaluations of the function and its first $p$ derivatives to compute such a
  point.
\end{itemize}

\noindent
{\bf Outline.} The paper is organized as
follows. Section~\ref{section:grad-Holder} considers the minimization of
regularized polynomials in Banach spaces. Section~\ref{algo-s} then introduces
the class of Banach spaces of interest and details our general adaptive
regularization algorithm for first-order minimization in these spaces,
while Section~\ref{complexity-s} analyzes its evaluation complexity. We conclude
the paper in Section~\ref{concl-s} with a brief discussion of the new results and perspectives.

\vspace*{1mm}
\noindent
\textbf{Notation} Throughout the paper, $\np{.}$ denotes the norm over the
space $\calV$. $\mathcal{B}(x,B)$ denotes the open ball centered at $x$ of
radius $B$.  $\mathcal{L} (\calV^{\otimes m} ; \Re) $ denotes the space of
multilinear continous functions from $\calV \times \calV \dots \times \calV$
to $\Re$ and  $\mathcal{L}^m_{sym} (\calV^{\otimes m} ; \Re)$ the subspace of
$\mathcal{L}^m (\calV^{\otimes  m} ; \Re) $ that is $m$-linear symmetric. For a
function $f$ defined from $\calV$ to $\Re$ that is $p$ times Fr\'echet
differentiable, $\nabla_x^k f(x) \in \mathcal{L}^k_{sym} (\calV^{\otimes k} ;
\Re)$ denotes the $k$-th  derivative tensor for $k \in \ii{p}$. $\nabla_x^1
f$ is an element of the dual space of $\calV$ denoted $\calV^\prime$.  The
symbol $\pd{\cdot,\cdot}$ denotes the dual pairing between $\calV$ and $\calV'$,
that is $\pd{y,x} \eqdef y(x)$, for $y\in \calV'$ and $x\in \calV$.
The norm in the dual space $\calV^\prime$ will be denoted as $\nd{.}$.
For $S \in \mathcal{L}_{sym}^m $, $S[v_1,v_2 \dots ,v_m] \in \Re$ denotes the
result of applying $S$ to the vectors $v_1, \dots ,v_m$. $S[v]^m$ is the
result of applying S to $m$ copies of the vector $v$ and $S[v]^l \in
\mathcal{L}^{m-l}_{sym} (\calV^{\otimes m-l} ; \Re) $ the result of applying
$l$ times the vector $v$. We define the norm in $\mathcal{L}^m_{sym}
(\calV^{\otimes  m} ; \Re)$ as
\beqn{Tensornormdef}
\|S\| \eqdef\sup_{\np{v_1}= \dots = \np{v_m} = 1}  | S[v_1, \dots, v_m ] |.
\eeqn
 
\numsection{Gradient descent with a H\"older regularization}\label{section:grad-Holder}

We start by considering the minimization, for $x$ in the Banach space $\calV$,
of the regularized polynomial functional of the form
\beqn{phi}
\phi(x) \eqdef \phi_0  + \sum_{\ell=1}^p \frac{1}{\ell!} S_\ell[x]^\ell   + h(x),
\eeqn
where $S_\ell \in \mathcal{L}^\ell_{sym} (\calV^{\otimes \ell})$ for $\ell \in
\ii{p}$ and $h$ is a general regularization term. Note that the sum of the two first
 terms of the right-hand side have the form of a Taylor expansion (in the
Fr\'{e}chet sense). The functions $\phi$ and  $h$ and the space $\calV$ 
are assumed to satisfy the following assumptions.
 
 \noindent
 \textbf{AS.1}
 \begin{itemize}
 \item[(i)]
   There exists $\phi_{\min} \in \Re$ such that, for all $x \in \calV$,
   $\phi(x) \geq \phi_{\min}$.  Moreover the set $\calD \eqdef \{ x \in
          \calV \, , \, \phi(x) \leq \phi(0)\}$ is bounded in the sense that
          $\sup_{x\in \mathcal{D}} \np{x} \leq \omega$ for some $\omega <
          \infty.$
          
\item[(ii)] $h$ is a convex differentiable function whose gradient satisfies the
local H\"older condition
\[ 
\forall \delta > 0, \, \forall x \in \mathcal{B}(0,\delta), \, \forall y \in \calV, \,
\nd{\nabla_x^1 h(x) - \nabla_x^1 h(y)}
\leq \sum_{i = 1}^k L_{i,\delta} \np{x-y}^{\beta_i-1}, 
\]
where, for $i\in \ii{k}$, $\beta_i >1$  and $L_{i,\delta}$
are a positive constants, the latter depending on $\delta$.
Moreover, $\beta_i \le 2$ for at least one $i\in\ii{k}$.
\item[(iii)] the space $\calV$ is reflexive.
\end{itemize}

\noindent
Observe that the condition stated in \textbf{AS.1}(ii) reduces to the standard
$\beta_1-1$-H\"{o}lder continuity of the gradients of $h$ whenever
$k=1$. Also note that, if all $\beta_i$ were strictly larger than two, $h$ would be affine.

We now use the property that $ S_\ell \in \mathcal{L}^\ell_{sym} (\calV^{\otimes \ell} ;
\Re)$  to derive an upper bound of $\phi(x+s)$ for all $x \in \mathcal{D}, s
\in \calV$. We then choose a specific $s$ to obtain the next result.

\llem{phi-upper}{There exists an integer $m \geq p$ and constants
  $\kappa_{i,\omega}>0$ ($i\in \ii{m}$) such that, for all $x \in \calD$, there exists a vector $d$ in $\calV$, 
 \beqn{graddsecant}
 \phi(x-td) \leq \phi(x) - \nd{\nabla_x^1 \phi(x)} t + \sum_{i=1}^m \kappa_{i,\omega} t^{\gamma_i},
 \eeqn
 where $1< \gamma_1 \leq \gamma_2 \leq \ldots \leq \gamma_m$ and $t \in \Re$.
}

\proof{
Successively using the binomial expansion, the convexity of $h$,
\req{Tensornormdef}, the fact that $\np{x}\leq \omega$ because $x\in\calD$ and
\textbf {AS.1}(ii), we derive that
\begin{align}
 \phi(x+s) &= \phi_0 + \sum_{\ell=1}^p \frac{1}{\ell!} S_l[x+s]^\ell + h(x+s), \nonumber \\
&= \phi(x)
 + \sum_{\ell=1}^p \frac{1}{(\ell-1)!} \pd{S_\ell[x]^{\ell-1} , s}
 + \sum_{\ell=2}^p \sum_{i=0}^{\ell-2} \frac{1}{\ell!}{\ell \choose i} S_\ell[x]^{i}[s]^{\ell-i} + h(x+s)-h(x), \nonumber \\
& \leq \phi(x) + \sum_{\ell=1}^p \frac{1}{(\ell-1)!} \pd{S_\ell[x]^{\ell-1} , s}
 + \sum_{\ell=2}^p \sum_{i=0}^{\ell-2} \frac{1}{\ell!}{\ell \choose i}\|S_\ell\| \np{x}^i \np{s}^{\ell-i}
 + \pd{\nabla_x^1 h(x+s) , s}, \nonumber \\
&\leq \phi(x) + \sum_{\ell=1}^p \frac{1}{(\ell-1)!} \pd{S_\ell[x]^{\ell-1} , s}
 + \sum_{\ell=2}^p \sum_{i=0}^{\ell-2} \frac{1}{\ell!}{\ell \choose i}\|S_\ell\| \omega^i \np{s}^{\ell-i}
 + \pd{\nabla_x^1 h(x+s), s}, \nonumber  \\ 
&\leq \phi(x) + \pd{\nabla_x^1 \phi(x) , s} + \sum_{\ell=2}^p \kappa_{\ell,\omega} \np{s}^\ell
 + \pd{\nabla_x^1 h(x+s) - \nabla_x^1 h(x) , s}, \nonumber   \\
&\leq \phi(x) + \pd{\nabla_x^1 \phi(x) , s} +\sum_{\ell=2}^p \kappa_{\ell,\omega} \np{s}^\ell +
 \sum_{\ell=1}^k L_{\ell,\omega} \np{s}^{\beta_\ell} \nonumber. 
\end{align}
Rearranging the last equation, we obtain that
\[
 \phi(x+s) \leq \phi(x) + \pd{\nabla_x^1 \phi(x), s} + \sum_{i=1}^m \kappa_{i,\omega} \np{s}^{\gamma_i},
\]
where the exponents $\gamma_i$ are in ascending order and strictly larger than
one. We now use the  reflexivity of $\calV$ to choose a $d \in \calV$ that
verifies both  $\pd{\nabla_x^1  \phi(x),d} = \nd{\nabla_x^1 \phi(x)} $ and
$\np{d} = 1$, we choose $s = -td$ in the last inequality so that \req{graddsecant} follows.
} 

\noindent
Looking at the steepest descent direction for minimizing \req{phi}, we are now
lead to consider \req{graddsecant} and to charaterize the minima of functions
of the form
\beqn{polynomialmodel}
 \Psi(t) \eqdef -\alpha t + \sum_{i = 1}^m \kappa_i t^{\gamma_i},
 \eeqn
for $t \in \Re_+$, $\alpha > 0$, $\kappa_i > 0$ and $1 < \gamma_1 \leq \gamma_2 \leq
\dots \leq \gamma_m$. This is the object of the next lemma.
 
 \llem{functionPsi}{
 	A function of the form \eqref{polynomialmodel} admits a unique minimum
        $t^\star$ over $\Re_+$ and
 	\beqn{descarg}
 	\Psi(t^\star) \leq -\min(\kappa_A \alpha^{\frac{\gamma_1}{\gamma_1 - 1 }} , \kappa_B \alpha^{\frac{\gamma_m}{\gamma_m - 1}}),
 	\eeqn
        where $\kappa_A$ and $\kappa_B$ depend on $\{\kappa_i\}_{i=1}^m$.
 }
 \proof{Let us consider $\Psi$ of the form \eqref{polynomialmodel}. Clearly,
   $\Psi$ is a strictly convex function as a sum of a linear function and a
   positive linear combination of powers strictly exceeding one. In
   addition, $\Psi^\prime(0) < 0$  and $\Psi^\prime(t) > 0$ for $t\in\Re_+$
   sufficiently large. Thus, a unique positive minimizer $t^\star$ exists such that
   $\Psi^\prime(t^\star) = 0$. Suppose first that $t^\star \geq 1$. Our
   problem then reduces to the minimization of $\Psi$ for $t \geq 1$. 
   Define
   \beqn{prop1}
   t_1
   \eqdef \left( \frac{\alpha}{\sum_{i=1}^m \kappa_i
     \gamma_i}\right)^{\sfrac{1}{\gamma_m - 1}}
   >0.
   \eeqn
   Because $t^\star \geq 1$ and $\psi^\prime$ is a non decreasing function,
   we obtain that
   $
   \psi^\prime(1) \leq  \psi^\prime(t^\star) = 0
   $
   and thus that
   \[
   - \alpha + \sum_{i=1}^m \kappa_i \gamma_i \leq 0,
   \]
   which, together with the definition of $t_1$ in \req{prop1}, implies that
   $
   (t_1)^{\gamma_m - 1 }\geq 1,
   $
   and the inequality $\gamma_m > 1$ then ensures $t_1\geq 1$.
   Using now the assumption that $1<\gamma_1 \leq \gamma_2 \leq \dots \leq
   \gamma_m$, we deduce that,  for $t \geq 1$,
   \[
   \Psi(t) \leq -\alpha t + \sum_{i=1}^m \kappa_i t^{\gamma_m},
   \]
   and thus, since $t_1\geq 1$,
   \begin{align*}
   \Psi(t_1) 
   &\leq - \alpha^{\sfrac{\gamma_m}{\gamma_m -1}} \left(\bigsum_{i=1}^m\kappa_i \gamma_i \right)^{\sfrac{-1}{\gamma_m-1}}
   + \alpha^{\sfrac{\gamma_m}{\gamma_m -1}}\bigfrac{\bigsum_{i=1}^m\kappa_i}
      {\left(\bigsum_{i=1}^m \kappa_i \gamma_i \right)^{\sfrac{\gamma_m}{\gamma_m - 1}}}, \\ 
   &\leq  -\alpha^{\sfrac{\gamma_m}{\gamma_m -1}}
       \left(\bigsum_{i=1}^m\kappa_i \gamma_i \right)^{\sfrac{-1}{\gamma_m-1}}
       \left(1 -\frac{\bigsum_{i=1}^m \kappa_i}{\bigsum_{i=1}^m \kappa_i \gamma_i }\right),\\
   & \eqdef - \alpha^{\sfrac{\gamma_m}{\gamma_m -1}} \, \kappa_A.
   \end{align*}
   As $\kappa_i > 0$ and $\gamma_i > 1$ for all $i$, we obtain that
   $\kappa_A > 0$ and hence $\Psi(t^\star) \leq \Psi(t_1) \leq - \kappa_A
   \alpha^{\sfrac{\gamma_m}{\gamma_m - 1}}$, which corresponds to the
   first term in the minimum of \eqref{descarg}.
 	
   Suppose now that $t^\star \leq 1$ and define
   \beqn{prop2}
   t_2 \eqdef \left( \frac{\alpha}{\sum_{i=1}^m \kappa_i \gamma_i}\right)^{\sfrac{1}{\gamma_1 - 1}}
   >0.
   \eeqn 
   As $\psi^\prime(t^\star) = 0$, $t^\star \leq 1$ and $1 < \gamma_1 \leq
   \gamma_2 \leq \dots \leq \gamma_m $,
   \[
   - \alpha + \sum_{i=1}^m \kappa_i \gamma_i (t^\star)^{\gamma_i - 1} = 0,
   \]
   and thus
   \[
   -\alpha  + \sum_{i=1}^m \kappa_i \gamma_i (t^\star)^{\gamma_1 - 1} \geq 0,
   \]
   which, with the definition of $t_2$ in \req{prop2} gives that 
   \[
 	(t^\star)^{\gamma_1 - 1} \geq t_2^{\gamma_1 - 1}. 	
   \]
   The inequality $\gamma_1 > 1$ then ensures that $t_2 \leq t^\star \leq 1$.   
   Using an argument similar to that used above but now for the case $t^\star
   \leq 1$, we deduce that, for all $t \leq 1$,
   \[
   \Psi(t) \leq -\alpha t + \sum_{i=1}^m \kappa_i t^{\gamma_1},
   \]
   and therefore, since $t_2\leq 1$,
   \begin{align*}
     \Psi(t_2)
     &\leq - \alpha^{\sfrac{\gamma_1}{\gamma_1 -1}} \left(\sum_{i=1}^m \kappa_i\gamma_i \right)^{\sfrac{-1}{\gamma_1 -1}}
     + \alpha^{\sfrac{\gamma_1}{\gamma_1 -1}} \bigfrac{\sum_{i=1}^m \kappa_i}
     {\left(\bigsum_{i=1}^m \kappa_i \gamma_i \right)^{\sfrac{\gamma_1}{\gamma_1 - 1}}}, \\ 
     &\leq  - \alpha^{\sfrac{\gamma_1}{\gamma_1 -1}}
     \left(\bigsum_{i=1}^m \kappa_i \gamma_i \right)^{\sfrac{-1}{\gamma_1 -1}} \left(1 -
     \bigfrac{\bigsum_{i=1}^m \kappa_i}{\bigsum_{i=1}^m \kappa_i \gamma_i  } \right),\\
     &\eqdef - \alpha^{\sfrac{\gamma_1}{\gamma_1 -1}} \kappa_B.
   \end{align*}
   Rewriting the last line gives that $\Psi(t^\star) \leq \Psi(t_2) \leq -
   \kappa_B \alpha^{\sfrac{\gamma_1}{\gamma_1 - 1}}$, which completes the proof.
 }

 \noindent
 This result suggest the following algorithm for minimizing functions of the
 form \req{phi}. 

 \algo{gradHolder}{A First-Order Gradient Algorithm for
                  Minimizing \\\hspace*{30mm}Regularized Polynomials}
{
        \begin{description}
        \item[Step 0: Initialization. ] Set $x_0= 0$ and $k= 0$.
        \item[Step 1: Compute a search direction.] Compute
          $\nabla_x^1\phi(x_k)\in \calV'$ and select a direction $d_k$ such that
          $\nd{\nabla_x^1\phi(x_k)} = \pd{\nabla_x^1\phi(x_k),d_k}$ and
          $\|d_k\|=1$. If $\nd{\nabla_x^1\phi(x_k)} = 0$, stop and return the
          sequence $(x_0 , x_1 \dots , x_k)$.
        \item[Step 2: Stepsize definition.] Compute $t_k$ a global
           minimizer of $\phi(x_k-td_k)$.
         \item[Step 3: Define the next iterate. ] Set $x_{k+1}=
           x_k-t_k d_k$, increment $k$ by one and return to Step~1.
        \end{description}
}

\noindent
Note that the selection of $d_k$ in Step~1 is possible because $\calV$ is
reflexive, and that the minimization in Step~2 is possible because it occurs in a
one-dimensional space.

We now prove the following convergence result.
      
\lthm{gradtozero}{
Suppose that $\phi$, $h$ and $\calV$ verify \textbf{AS.1} and let
$\{x_k\}_{k\geq 0}$ be the sequence generated by
Algorithm~\ref{gradHolder}. Then
\[
\phi(x_{k+1}) < \phi(x_k) \tim{ for all } k\geq 0
\]
and  either the algorithm terminates in a finite number of iterations
  with an iterate $x_k$ such that $\nabla_x^1 \phi(x_k)=0$, or
\[
\lim_{k\rightarrow \infty} \nd{\nabla_x^1 \phi(x_k)} = 0.
\]
}

\proof{
Recall that $ \mathcal{D} = \{ x \in \calV \, \mid \, \phi(x) \leq\phi(x_0)=\phi(0)\}$
and that inequality \eqref{graddsecant} is valid if $x \in \mathcal{D}$.
Since the left hand-side of the inequality \eqref{graddsecant} for $x=x_0$ verifies the
conditions of Lemma~\ref{functionPsi}, and  denote by $t_0^\star$ the
minimizer of Lemma~\ref{functionPsi}. We may apply this lemma and deduce that,
\[ 
\phi(x_1 ) \leq \phi(x_0-t_0^\star d_0 ) \leq \phi(0) -  \min (\kappa_A \nd{\nabla_x^1
  \phi(0)}^{\sfrac{\gamma_1}{\gamma_1 -1}} , \kappa_B \nd{\nabla_x^1
  \phi(0)}^{\sfrac{\gamma_m}{\gamma_m -1}}),
\]
where now $\kappa_A$ and $\kappa_B$ are strictly positive and depend on
$\omega$ (the radius of $\calD$) and the Lipschitz constant $L_{i,\omega}$,
themselves depending on $\omega$. As $\nd{\nabla_x^1\phi(x_0)} > 0 $,
$\phi(x_1) < \phi(x_0)= \phi(0)$ and therefore $x_1 \in \calD$.

Suppose now that $x_{k-1}\in\calD$ and that $\nd{\nabla_x^1 \phi(x_{k-1}) } >0$.  We may again apply
Lemma~\ref{functionPsi} to the left handside of inequality \eqref{graddsecant}
with $x$ chosen as $x_{k-1}$ and by denoting $t_{k-1}^\star$ the minimizer of
the left hand-side, we deduce that
\[ 
\phi(x_k) \leq \phi(x_{k-1}-t_{k-1}^\star d_{k-1} ) \leq \phi(x_{k-1}) -  \min (\kappa_A \nd{\nabla_x^1
  \phi(x_{k-1})}^{\sfrac{\gamma_1}{\gamma_1 -1}} , \kappa_B \nd{\nabla_x^1
  \phi(x_{k-1})}^{\sfrac{\gamma_m}{\gamma_m -1}}),
\]
thus $x_k$ and the complete sequence $\{x_k\}_{k\geq 0}$ belong to
$\calD$ and the first conclusion of the theorem holds. To prove the
second part,
we first note that the definition of the algorithm ensures the identity
$\nabla_x^1\phi(x_k)=0$ whenever termination occurs after a finite number of
iterations. Assume therefore that the algorithm generates an 
infinite sequence of iterates and that
\beqn{absurd}
\nd{\nabla_x^1 \phi(x_{k_i})} \geq \epsilon,
\eeqn
for some $\epsilon> 0$ and some subsequence $\{k_i\}_{i=1}^\infty$.  Summing over all
iterations ${k_i}$ and using \textbf{AS.1}(i), we obtain that  
\begin{align*}
+\infty > \phi(0)-\phi_{\min}
&\geq \sum_i \min (\kappa_A \nd{\nabla_x^1\phi(x_{k_i})}^{\sfrac{\gamma_1}{\gamma_1 -1}} ,
  \kappa_B \nd{\nabla_x^1 \phi(x_{k_i})}^{\sfrac{\gamma_m}{\gamma_m -1}}),\\
&\geq  \sum_{i} \min[ \kappa_A\epsilon^{\sfrac{\gamma_1}{\gamma_1 -1}},
                   \kappa_B\epsilon^{\sfrac{\gamma_m}{\gamma_m-1}} ],
\end{align*}
which is a contradiction since the right-hand side diverges to
$+\infty$. Hence \req{absurd} cannot hold and the second conclusion of the
theorem is valid.
} 

Thus a vanilla gradient-descent algorithm applied to a $p$-th
degree polynomial augmented by a convex regularization term with H\"older
gradient will yield asymptotic first-order stationarity.

\numsection{An adaptive regularization algorithm in Banach spaces}~\label{algo-s}

\noindent
We now consider developing an adpative regularization method for finding
first-order points for the problem
\beqn{problem}
\min_{x\in\calV}f(x),
\eeqn
and make our assumptions on the problem more precise.

\noindent
\textbf{AS.2} $f$ is $p$ times continuously Fr\'echet differentiable with $p \geq 1$. 

\noindent
{\bf AS.3} There exists a constant $f_{\rm low}$ such that
$f(x) \geq f_{\rm low}$ for all $x\in\calV$.

\noindent
\textbf{AS.4} The $p$-th derivative tensor $\nabla_x^p f(x) \in \mathcal{L}
(\calV^p; \Re) $ is globally H\"older continuous, that is, there exist
constants $L > 0$ and $\beta \in (0,1]$ such that  
\beqn{holderpth}
\|\nabla_x^p f(x) - \nabla_x^p f(y)\| \leq L \np{x-y}^\beta, \, \text{ for all } x,y \in \calV .
\eeqn

\noindent
For brevity, \textbf{AS.2} and \textbf{AS.4} will be denoted by $f \in \mathcal{C}^{p,\beta}(\calV;\Re)$.

\noindent
Let $T_{f,p}(x,s)$ be the Taylor series of the functional $f(x+s)$ truncated at order $p$. 
\beqn{taylor model}
T_{f,p}(x,s) \eqdef f(x) + \sum_{l=1}^p \frac{1}{l!} \nabla_x^l f(x) [s]^l.  
\eeqn
The gradient $\nabla_x^1 f(x)$  belongs to the dual space $\calV^\prime$ and
will be denoted by $g(x)$. Thus, for a requested accuracy $\epsilon \in
(0,1]$, we are interested in finding an $\epsilon$-approximate first-order
critical point, that is a point $x_\epsilon$ such that $\nd{g(x_\epsilon)}
\leq \epsilon$. 

\subsection{Smooth Banach spaces}

\noindent 
In a generic Banach space, we can only ensure
``a decrease principle'' as stated in \cite[Theorem 5.22]{Clar13}. To
obtain more conclusive results, we need to introduce additional
assumptions. We choose to work with the class of \textit{uniformly q
smooth Banach spaces}. For the sake of completeness, we briefly recall the context.
 
\noindent
Given a Banach space $\calV$, we first define its module of smoothness, for $t
\geq 0$, by
\beqn{unifdef}
\rho_{\calV}(t)
\eqdef \sup_{\np{x} = 1 \, , \np{y} = t} \left\{ \frac{\np{x+y} + \np{x-y}}{2} - 1 \right\},
\eeqn
and immediately deduce from the triangular inequality that $\rho_{\calV}(t)
\leq t$. We now say that $\calV$ is a uniformly smooth Banach space if and
only if $\lim_{t \to 0 } \frac{\rho_{\calV}(t)}{t} = 0$. 
Going one step further,  we say that a Banach space $\calV$ is
\textit{uniformly q smooth} for some $q \in (1,2]$ if and only if 
\beqn{qsmoothdef}
\exists \kappa_\calV >0,  \, \rho_{\calV}(t) \leq \kappa_\calV t^q.
\eeqn
\noindent
It is easy to see that, if $\calV$ is uniformly $q$ smooth, it is also
uniformly $q^\prime$ smooth for all $1 < q^\prime < q$. Indeed, one can easily
show\footnote{If $t\in[0,1]$ this follows from \req{qsmoothdef} and $q^\prime <
  q$. If $t > 1$, $\rho_{\calV}(t) \leq t \leq t^{q^\prime}$.} that
$\rho_{\calV}(t) \leq \max(1,\kappa_\calV) t^{q^\prime}$ from definition
\eqref{unifdef} and inequality \eqref{qsmoothdef}. 

We motivate our choice of this particular class of Banach spaces by giving a
few examples. $L^p(\Re), \, 1 < p < \infty,$ are uniformly smooth Banach
spaces. In particular, $L^p(\Re)$ is uniformly 2 smooth for $p \geq 2$ and
uniformly $p$ smooth for $1< p \leq 2$. The same results apply for $\ell^p$
and the Sobolev spaces $W_m^p(\Re)$ \cite{XuRoac91}. Moreover, all Hilbert
spaces are 2 smooth Banach.

\llem{2smooth}{
  Let $\calH$ be a Hilbert space. Then $\calH$ is a 2 smooth Banach space with
	\beqn{2banach}
	\rho_{\calH}(t) \leq \frac{t^2}{2}.
	\eeqn}  
\proof{
Because of the definition of $\rho_{\calH}$ in \req{unifdef}, we have that
\begin{align*}
\rho_{\calH}(t) &= \sup \left\{ \frac{\np{x+y}+\np{x-y}}{2} - 1 , \np{x}=1 \, \, , \np{y}=t \right\},   \\
&=  \sup \left\{ \frac{\sqrt{1+t^2 +2 \pd{y,x}} + \sqrt{1+t^2 -2 \pd{y,x}} }{2} -1 \, , \np{x}=1 \, \, , \np{y}=t\right\}. 
\end{align*}
Thus, when maximizing over $\pd{y,x} \in [-t, t]$,
\[
\rho_{\calH}(t) = \sqrt{1+t^2} -1 = \frac{t^2}{\sqrt{1+t^2} +1} \leq \frac{t^2}{2}.
\]
} 

\noindent
One might wonder if it is possible for the $q$ smooth order to be strictly
superior to $2$ in \eqref{qsmoothdef}. We now show that this is impossible. 
Indeed, for any Banach space $\calV$, we have that, $\rho_{\calV}(t) \geq
\rho_{\calH}(t) = \frac{t^2}{\sqrt{1+t^2}+1}$ \cite{XuRoac91}. Suppose now  
$\rho_{\calV}(t) \leq c t^m$ with $m > 2$. Using the last two inequalities, we
obtain that: $c t^{m-2} \geq \frac{1}{\sqrt{1+t^2}+1}$ for all $t$ strictly
positive. But this  inequality is impossible for small enough $t$ and hence our
supposition about $ m$ is false and $m \in (1,2]$.

From here on, we assume that 

\vspace*{1mm}
\noindent
\textbf{AS.5} $\calV$ is a uniformly $q$ smooth space.
\vspace*{1mm}

\noindent
Uniformly smooth Banach spaces are also reflexive (See
\cite[Proposition 1.e.3, p61]{XuRoac91}), so that \textbf{AS.1}(iii) automatically holds. 
Let us now define the set
\beqn{sousdiffp}
J_p(x)
\eqdef \left\{ v^{*} \in \calV^{*} \, , \,\pd{v^*,x} = \np{x}^p \, , \,\nd{v^{*}} = \np{x}^{p-1} \right\}.
\eeqn
\noindent
It is known \cite{Xu91} that $J_p(x)$ is the subdifferential of the functional
$\frac{1}{p} \np{ \,\cdot \,}^p \, , \, p\geq 1$ at $x$.

\noindent
We may now introduce another characterization of uniform smoothness.

\lthm{unifromsmoothcara}{
Let 
\[
\mathcal{F} \eqdef \{\psi:\Re \to \Re \mid \psi(0) = 0, \psi \text{ is
  convex, non decreasing and } \exists \kappa_\calF > 0 \mid \psi(t) \leq \kappa_\calF \rho_{\calV}(t)\}.
\]
Then, for any $1 < p < \infty$, the following statements are equivalent. 
\begin{itemize}
\item[(i)] $\calV$ is a uniformly smooth Banach space. 
\item [(ii)] $J_p$ is single valued and there exists
   $\varphi_p(t) = \frac{\psi_p(t)}{t}$ where $\psi_p \in \calF$ and such that
   \beqn{subdiffcontrol}
   \nd{J_p(x) - J_p(y)}
   \leq \max(\np{x} , \np{y})^{p-1}  \varphi_p\left(\frac{\np{x-y}}{\max(\np{x}, \np{y})}\right).
   \eeqn
\end{itemize}}
\proof{\cite[Theorem~2]{XuRoac91}.}
As we will be only working with $\np{.}^p$ for $p > 1$ in the rest of the
paper, we define $J_p(x)$ as the unique value in the set \eqref{sousdiffp}. As
the subdifferential of $\np{.}^p$ reduces to a singleton for $p > 1$ and
$\np{.}^p$ is a convex function, $\np{.}^p$ is Fr\'echet differentiable for $p
> 1$ since it verifies \cite[Condition~4.16]{Clar13}. The reader
is referred to \cite{Xu91} or \cite{XuRoac91} for more extensive coverage of
characterizations of the norm in uniformly smooth Banach spaces.

For all $\ell > 1$, we now prove an upper bound of the norm of
$\nd{J_\ell(x)- J_\ell(y)}$ in terms of $\np{x-y}$ in a uniform $q$ smooth
Banach space. Let us first remind the useful inequality
$(x+y)^r \leq \max(1,2^{r-1}) ( x^r + y^r)$ for all $x,y \geq 0$ and all
$r \ge 0$, before stating the next crucial lemma.  

\llem{Ball Majorisation}{Suppose that $\calV$ is a uniformly $q$ smooth
  Banach space and that $x\in \mathcal{B}(0,\omega)$.
  Then for all $\ell > 1$, there exist constants $\kappa_\omega, \kappa_\ell > 0$ such
  that
  \beqn{Subgradcontrol}
	\nd{J_{\ell}(x) - J_{\ell}(y) } \leq \kappa_\omega \np{x-y}^{\min[q,\ell]-1} + \kappa_\ell \np{x-y}^{\ell-1},
  \eeqn
  where $\kappa_\omega$ and $\kappa_\ell$ depend only on $\omega$, $\ell$,
  $\kappa_\calF$ and $\kappa_\calV$.
}
\proof{
	As $\ell > 1$, if $q > \ell$, we can use our remark above and decrease the
        $q$ smooth order until $q' = \min[q,\ell] \leq \ell$.
        We now develop the upper bound (ii) of Theorem \ref{unifromsmoothcara}
        and use the definition of the set $\mathcal{F}$ to derive that
	\begin{align*}
	\nd{J_{l}(x) - J_{l}(y) } &\leq \max(\np{x} , \np{y})^{\ell-1}
        \kappa_\calF \kappa_\calV \left(\frac{\np{x-y}}{\max(\np{x}, \np{y})}\right)^{q'-1}, \\ 
	&\leq \max(\np{x} , \np{y})^{l-q'} \kappa_\calF \kappa_\calV \np{x-y}^{q'-1}.
	\end{align*}
	Using now the inequalities $\max(\np{x} , \np{y}) \leq  \np{x} +\np{x-y} 
        $ and $\ell \geq q'$, 
        we obtain that
	\begin{align*}
	\nd{J_{l}(x) - J_{l}(y) } &\leq \kappa_\calF \kappa_\calV ( \np{x} + \np{x-y} )^{\ell-q'} \np{x-y}^{q'-1}, \\
	& \leq \kappa_\calF \kappa_\calV \max(1,2^{\ell-q'-1})
        (\np{x}^{\ell-q'} + \np{x-y}^{\ell-q'}) \np{x-y}^{q'-1}, \\
	&\leq \kappa_\calF \kappa_\calV \max(1,2^{\ell-q'-1})\, \omega^{\ell-q'} \np{x-y}^{q'-1}\\
        & \hspace*{20mm}+ \kappa_\calF \kappa_\calV \max(1,2^{\ell-q'-1})  \np{x-y}^{\ell-1}, \\
	&\leq \kappa_\omega \np{x-y}^{q'-1} + \kappa_\ell \np{x-y}^{\ell-1}.    
	\end{align*}
} 

\subsection{The \al{AR$p$-BS} algorithm}

Adaptive regularization methods are iterative schemes which compute a step form
an iterate $x_k$ by building, for $f\in\mathcal{C}^{p,\beta}(\calV;\Re)$, a
regularized model $m_k(s)$ of $f(x_k+s)$ of the form
\beqn{model}
m_k(s) \eqdef T_{f,p}(x_k,s) + \frac{\sigma_k}{(p+\beta)!}  \np{s}^{p+\beta} \, , \, p\geq 1.
\eeqn
As in \cite{CartGoulToin18a} but at variance with \cite{CartGoulToin16},
we will assume here that $\beta$, the degree of H\"older continuity of the
$p$-th derivative tensor of $f$, is known. The $p$-th order Taylor series is
``regularized'' by adding the term $\frac{\sigma_k}{(p+\beta)!}
\np{s}^{p+\beta}$, where $\sigma_k$ is known as the ``regularization
parameter''. This term guarantees that the functionnal $m_k(s)$ is bounded below and thus
makes the procedure of finding a step $s_k$ by (approximately) minimizing
$m_k(s)$ well-defined. In our uniform $q$ smooth setting, $m_k(s)$ is
Fr\'echet differentiable but this is unfortunately insufficient to derive
results on the Lipschitz continuity of its gradient, which makes the use of
more standard gradient-descent methods impossible.

\noindent
Our proposed algorithm is similar in spirit to \al{ARC} \cite{CartGoulToin11d}
and proceeds as follows. At a given iterate $x_k$, a step $s_k$ is first
computed by approximately minimizing \req{model}. Once the step is computed,
the value of the objective functional at the trial point $x_k+s_k$ is then
evaluated.  If the decrease in $f$ from $x_k$ to $x_k+s_k$ is comparable to
that predicted by the $p$-th order Taylor series, the trial point is accepted
as the new iterate and the regularization parameter is (possibly) reduced. If
this is not the case, the trial point is rejected and the regularization
parameter is increased.  The resulting algorithm is formally stated as the
\al{AR$p$-BS} algorithm on the next page.

 \algo{ARpBS}{$p$-th order adaptive regularization in a uniform $q$ smooth Banach Space (\tal{AR$p$-BS})}{
	\begin{description}
		\item[Step 0: Initialization: ] An initial point $x_0\in \calV$, a regularization
		parameter $\sigma_0$ and a requested final gradient accuracy
		$\epsilon \in (0,1]$ are given. The constants
		$\eta_1$, $\eta_2$, $\gamma_1$, $\gamma_2$, $\gamma_3$, $\chi \in (0,1) $, and $\sigma_{\min}$
		are also given such that
		\beqn{banach-hyparam}
		\sigma_{\min} \in (0, \sigma_0], 
		0 < \eta_1 \leq \eta_2 < 1 
		\tim{ and } 0< \gamma_1 < 1 < \gamma_2 < \gamma_3.
		\eeqn
		Compute $f(x_0)$ and set $k=0$.
		\item[Step 1: Check for termination: ] Terminate with $x_\epsilon = x_k$ if
		\beqn{stopcondhilb}
		\nd{g(x_k)} \leq \epsilon.
		\eeqn
		\item[Step 2: Step calculation: ] Compute a step $s_k$  which
		sufficiently reduces the model $m_k$ in the sense that
		\beqn{qsmoothdescent}
		m_k(s_k) < m_k(0),
		\eeqn
		and
		\beqn{qsmoothstep}
                \nd{\nabla_s^1 m_k(s_k)} \leq \max\left[\chi  \epsilon, \,\theta \np{s_k}^{p+\beta-1}\right].
		\eeqn
		\item[Step 3: Acceptance of the trial point. ]
		Compute $f(x_k+s_k)$ and define 
		\beqn{qsmoothrhok-def}
		\rho_k = \frac{f(x_k) - f(x_k+s_k)}{T_{f,p}(x_k,0)-T_{f,p}(x_k,s_k)}.
		\eeqn
		If $\rho_k \geq \eta_1$, then define
		$x_{k+1} = x_k + s_k$; otherwise define $x_{k+1} = x_k$.
		\item[Step 4: Regularization parameter update. ]
		Set
		\beqn{qsmoothsigma-update}
		\sigma_{k+1} \in \left\{ \begin{array}{ll}
			{}[\max(\sigma_{\min}, \gamma_1\sigma_k), \sigma_k ]  & \tim{if} \rho_k \geq \eta_2, \\
			{}[\sigma_k, \gamma_2 \sigma_k ]          &\tim{if} \rho_k \in [\eta_1,\eta_2),\\
			{}[\gamma_2 \sigma_k, \gamma_3 \sigma_k ] & \tim{if} \rho_k < \eta_1.
		\end{array} \right.
		\eeqn
		Increment $k$ by one and go to Step~1.
	\end{description}
}

\noindent
While the \al{AR$p$-BS} algorithm follows the main lines of existing
\al{AR$p$} methods \cite{CartGoulToin11d,BirgGardMartSantToin17}.
Because we are in an infinite dimensional space, the
existence of a minimizer of $m_k(s)$ may not be guaranteed and hence a point
$s^\star$ such that $\nabla_s^1 m_k(s^\star) = 0$ may not exist. As a
consequence, standard proofs that a step satisfying both \req{qsmoothdescent}
and \req{qsmoothstep} exists no longer apply. We thus need to check that this is
still the case in our context. This is achieved using Algorithm~\ref{gradHolder}.

\lthm{step-exists}{
Suppose that \textbf{AS.2}, \textbf{AS.4} and \textbf{AS.5} hold.
Suppose also that $\nd{g(x_k)} > 0$.  Then a step satisfying both
\req{qsmoothdescent} and \req{qsmoothstep} always exists.
}

\proof{
  First note that \textbf{AS.2} and \textbf{AS.4} imply that $p+\beta > 1$.
In order to apply Algorithm~\ref{gradHolder} to the problem of
  minimizing \req{model}, we just need to prove that $m_k(s)$ satisfies \textbf{AS.1} of Section 2. 
We have that
\[
m_k(s)
\geq m_k(0) - \sum_{i=1}^p \|\nabla_x^i f(x)\| \np{s}^i
+ \frac{\sigma_k}{(p+\beta)!} \np{s}^{p+\beta} \to \infty
\text{ as } \np{s} \to \infty,
\]
\noindent 
and thus $m_k$ is a coercive functional  verifying  \textbf{AS.1}(i). Lemma
\ref{Ball Majorisation} (applied with $k=2$, $\delta=\omega$,
$\ell = p+\beta$, $L_{1,\delta}= \kappa_\ell$, $\beta_1= \min[q,\ell] \in (1,2]$,
$L_{2,\delta}= \kappa_\omega$ and $\beta_2=\ell>1$)
then ensures that $\np{.}^{p+\beta} $ satisfies
\textbf{AS.1}(ii). We already noted that, being uniformly smooth, $\calV$ must be
reflexive, which ensures that \textbf{AS.1}(iii) holds.  All the
requirements of \textbf{AS.1} in Section 2 are therefore met and, since $\nabla_s^1
m_k(0) = g(x_k)$, Theorem \ref{gradtozero} applies to the functional $m_k(s)$. As
a consequence, a suitable step $s_k$ such that $m_k(s_k) < m_k(0)$ and
$\nd{\nabla_s^1 m_k(s_k)} \leq \chi \epsilon $ exists.  
}

\noindent
Observe that equation \eqref{graddsecant} and the fact
that $\gamma_1 = \min[q , p+\beta]$ and $\gamma_m = p+\beta$ (all the other powers ranging
from $2$ to $p$), imply that, for our iterative gradient descent,
\[
\lim_{i \rightarrow \infty}
\min\left[\kappa_A \nd{\nabla_s^1 m(s_i)}^{\frac{\min[q , p+\beta]}{\min[q , p+\beta]-1}},
     \kappa_B\nd{\nabla_s^1 m(s_i)}^{\frac{p+\beta}{p+\beta-1}} \right] = 0.
\]
As a consequence, the first term in the minimum indicates that the smoother
the space, the faster the convergence for $p\geq 2$.

\noindent
Following well-established practice, we now define
\[
\calS \eqdef \{ k \geq 0 \mid x_{k+1} = x_k+s_k \} = \{ k \geq 0 \mid \rho_k \geq \eta_1 \},
\]
the set of indexes of ``successful iterations'', and
\[
\calS_k \eqdef \calS \cap \ii{k},
\]
the set of indexes of successful iterations up to iteration $k$.
We also recall a well-known result bounding the total number of iterations in
terms of the number of successful ones.

\llem{SvsU}{
Suppose that the \al{AR$p$-BS} algorithm is used and that $\sigma_k \leq
\sigma_{\max}$ for some $\sigma_{\max} >0$. Then
\beqn{unsucc-neg}
k \leq |\calS_k| \left(1+\frac{|\log\gamma_1|}{\log\gamma_2}\right)+
\frac{1}{\log\gamma_2}\log\left(\frac{\sigma_{\max}}{\sigma_0}\right).
\eeqn
}

\proof{See \cite[Theorem~2.4]{BirgGardMartSantToin17}.}

\numsection{Evaluation complexity for the \tal{AR$p$-BS} algorithm}\label{complexity-s}

\noindent
Before discussing our analysis of evaluation complexity, we first restate some
classical lemmas of \al{AR$p$} algorithms, starting with H\"older error bounds.

\llem{lipschitz}{
  Suppose that$f \in \mathcal{C}^{p,\beta}(\calV;\Re)$  holds and that $k\in \calS$. Then
  \beqn{Lip-f}
  |f(x_{k+1})- T_{f,p}(x_k,s_k)| \leq \frac{L}{(p+\beta)!}\np{s_k}^{p+\beta},
  \eeqn
  and
  \beqn{Lip-g}
  \nd{g_{k+1}-\nabla_s^1 T_{f,p}(x_k,s_k)}  \leq \frac{L}{(p-1+\beta)!} \np{s_k}^{p-1+\beta}.
  \eeqn
}

\proof{This is a direct extension of \cite[Lemma~2.1]{CartGoulToin20b}
  since the proof in this reference only involves \textbf{AS.2}, \textbf{AS.4} and unidimensional integrals. }

From now on, the analysis follows that presented in
\cite{BirgGardMartSantToin17} quite closely.
    
\llem{model-decrease}{
  \beqn{Tdecr}
  \Delta T_{f,p}(x_k,s_k) \eqdef T_{f,p}(x_k,0)-T_{f,p}(x_k,s_k)
  \geq \frac{\sigma_k}{(p+\beta)!} \np{s_k}^{p+\beta}.
  \eeqn
}

\proof{Direct from \req{qsmoothdescent} and \req{model}.}

\llem{sigma-max}{
  Suppose that $f \in \mathcal{C}^{p,\beta}(\calV;\Re)$.  Then, for all $k\geq 0$,
  \beqn{sigma-upper}
  \sigma_k\leq\sigma_{\max} \eqdef \gamma_3\max\left[\sigma_0,\frac{L}{(1-\eta_2)}\right].
  \eeqn
}

\proof{See \cite[Lemma~2.2]{BirgGardMartSantToin17}.
Using \req{qsmoothrhok-def}, \req{Lip-f},
and \req{Tdecr}, we obtain that
\[
|\rho_k-1|
\leq \bigfrac{(p+\beta)!|f(x_k+s_k)-T_{f,p}(x_k,s_k)|}{\sigma_k \np{s_k}^{p+\beta}}
\leq \bigfrac{L}{\sigma_k}.
\]
Thus, if $\sigma_k \geq L/(1-\eta_2)$, then $\rho_k \geq \eta_2$ ensures that iteration
$k$ is successful and \req{qsmoothsigma-update} implies that
$\sigma_{k+1}\leq \sigma_k$. The mechanism of the algorithm then guarantees that
\req{sigma-upper} holds.
} 

\noindent
The next lemma remains in the spirit of
\cite[Lemma~2.3]{BirgGardMartSantToin17}, but now takes the condition
\req{qsmoothstep} into account.

\llem{useful}{
  Suppose that $f \in \mathcal{C}^{p+\beta}(\calV; \Re)$ holds and
  that $k\in\calS$ before termination. Then
  \beqn{crucial}
  \np{s_k}^{p-1+\beta}
  \geq \epsilon \min\left[\frac{(1-\chi) (p+\beta-1)!}{L + \sigma_{\max}} ,
                 \frac{(p+\beta-1)! }{L+\sigma_{\max} + \theta (p+\beta-1)! }\right].
  \eeqn
}

\proof{
  Successively using the fact that termination does not occur at iteration
  $k$ and condition \req{qsmoothstep}, we deduce that
  \[
  \begin{array}{lcl}
    \epsilon & < &\nd{g(x_{k+1})},  \\*[2ex]
  & \leq &  \nd{g(x_{k+1})- \nabla_s^1 T_{f,p}(x_k,s_k)}
      + \nd{ \nabla_s^1 m_k(s_k)}
      +   \bigfrac{\sigma_k}{(p+\beta-1)!} \nd{J_{p+\beta}(s_k)}, \\*[2ex]
  & \leq & \bigfrac{L}{(p-\beta+1)!} \np{s_k}^{p-1+\beta}
      + \max\left[\chi \epsilon , \theta \np{s_k}^{p-\beta+1}\right]
      + \bigfrac{\sigma_k}{(p+\beta-1)!} \np{s_k}^{p+\beta-1}.
  \end{array}
  \]
   By treating each case in the maximum separately, we obtain that either
   \[
   (1-\chi) \epsilon \leq \left(\frac{L}{(p+\beta-1)!} + \frac{\sigma_k}{(p+\beta-1)!} \right) \np{s_k}^{p-1+\beta},
   \]
   or
   \[
   \epsilon \leq \left(\frac{L}{(p+\beta-1)!} + \frac{\sigma_k}{(p+\beta-1)!} + \theta \right) \np{s_k}^{p-1+\beta}. 
   \]
   Combining the two last inequalities gives that
   \[ 
   \np{s_k}^{p-1+\beta} \geq \min\left[\frac{(1-\chi)\epsilon (p+\beta-1)!}{L + \sigma_{\max}} , \frac{(p+\beta-1)! \epsilon}{L+\sigma_{\max} + \theta (p+\beta-1)! }\right].
   \]

  This in turn directly implies \req{crucial}.
} 

\noindent
We may now resort to the standard ``telescoping sum'' argument to obtain the
desired evaluation complexity result.

\lthm{complexity}{
  Suppose that \textbf{AS.2--AS.5} hold.   Then the \al{AR$p$-BS} algorithm requires at most
  \[
  \kappa_{\rm ARpBS}\,\frac{f(x_0)-f_{\rm low}}{\epsilon^{\frac{p+\beta}{p+\beta-1}}},
  \]
  successful iterations and evaluations of $\{\nabla_x^i f\}_{i=1,2,\dots,p}$ and at most
  \[
  \kappa_{\rm ARpBS}\,\frac{f(x_0)-f_{\rm low}}{\epsilon^{\frac{p+\beta}{p+\beta-1}}}
  \left(1+\frac{|\log\gamma_1|}{\log\gamma_2}\right)+
  \frac{1}{\log\gamma_2}\log\left(\frac{\sigma_{\max}}{\sigma_0}\right),
  \]
  evaluations of $f$ to produce a vector $x_\epsilon\in \calV$ such that
  $\nd{g(x_\epsilon)} \leq \epsilon$, where
  \[
  \kappa_{\rm ARpBS}
  =\frac{(p+\beta-1)!}{\eta_1\sigma_{\min} }
   \min\left[\frac{(1-\chi)\epsilon (p+\beta-1)!}{L + \sigma_{\max}} ,
             \frac{(p+\beta-1)! \epsilon}{L+\sigma_{\max} + (p+\beta-1)!\theta}
             \right]^{\frac{p+\beta}{p+\beta-1}}.
  \]
}

\proof{
  Let $k$ be the index of an iteration before termination.  Then, using \textbf{AS.3},
  the definition of successful iterations, \req{Tdecr} and \req{crucial},
  and the fact that computing an appropriate step is of constant order of
  complexity, we obtain that
  \[
  f(x_0)-f_{\rm low}
  \geq \bigsum_{i=0,i \in\calS}^k f(x_i)-f(x_{i+1})
  \geq \eta_1 \bigsum_{i \in\calS_k} \Delta T_{f,2}(x_i,s_i)
  \geq \bigfrac{|\calS_k|}{\kappa_{\rm ARpBS}} \epsilon^{\frac{p+\beta}{p+\beta-1}}.
  \]
  Thus
  \[
  |\calS_k| \leq \kappa_{\rm ARpBS}\,\frac{f(x_0)-f_{\rm low}}{\epsilon^{\frac{p+\beta}{p+\beta-1}}},
  \]
  for any $k$ before termination. The first conclusion follows since the
  derivatives are only evaluated once per successful iteration.  Applying now
  Lemma~\ref{SvsU} gives the second conclusion.

} 
\noindent
Theorem~\ref{complexity} extends the result of \cite{BirgGardMartSantToin17}
in the case $\beta = 1$ and some results of \cite{CartGoulToin20b} to
uniform $q$ smooth Banach spaces. We recall that  $L^p$, $\ell^p$ and $W_m^p$
are uniform $q$ smooth spaces for $1< p <\infty$, and hence that
Lemma~\ref{Ball Majorisation} and Theorem~\ref{complexity} apply in these spaces.
We may also consider the finite dimensional case where $\Re^n$ is equipped
with the norm $\| x \|_r = \left( \sum_{i=1}^n |x_i|^r
\right)^{\sfrac{1}{r}}$.  We know that, for all $1<r<\infty$, this is a
uniform $\min(r,2)$ smooth space, and therefore Theorem~3.5 again applies.
We could of course have obtained convergence of the adaptive regularization algorithm in
this case using results for the Euclidean norm and introducing
norm-equivalence constants in our proofs and final result, but this is avoided
by the approach presented here.  This could be significant when the dimension
is large and the norm-equivalence constants grow.

\numsection{Discussion}\label{concl-s}

We have proposed a generalized H\"older condition and a
gradient-descent algorithm for minimizing polynomial
functionals with a general convex regularization term in Banach spaces, and have
applied this result to show the existence of a suitable step in an adaptive
regularization method for unconstrained minimization in $q$ smooth Banach spaces.
We have also analyzed the evaluation complexity of this latter algorithm and
have shown that, under standard assumptions, it will find an $\epsilon$-approximate
first-order critical point in at most
$\calO\big(\epsilon^{-\frac{p+\beta}{p+\beta-1}}\big)$ evaluations of the functional and
its first $p$ derivatives, which is identical to the bound
known for minimization in  (finite-dimensional) Euclidean spaces. Since these
bounds are known to be sharp \cite{CartGoulToin18a}, so is ours.

It would be interesting to consider convergence to second-order points, but
the infinite dimensional framework causes more difficulties.
Indeed, considering second-order derivatives as in \cite{CartGoulToin12} is
impossible since we do not know if a power of the norm is twice
differentiable. As an example, consider $L^r([0,1])$ for $p > 1$, where
\[ 
\nabla_f^1 \left(\frac{\| f\|^p_{L^r([0,1])}}{p} \right) = \| f\|_{L^r([0,1])}^{p-r} f  | f|^{r-2}.   	 
\]
The right-hand side of the last equation involves an absolute value which is
only differentiable for specific values of $r$. It is interesting to
study the case of $r=2$ with the objective of extending our analysis to the second order.
Another line of future work is  to extend these results to metrizable spaces
(using the Bergman divergence or the Wasserstein distance) and to the complexity of second
order adaptive regularization in an infinite-dimensional Hilbert space.

{\footnotesize


}

\end{document}